\documentclass[12pt,A4]{article}
\usepackage{enumerate}
\usepackage{booktabs} 
\usepackage{psfrag}
\usepackage{graphicx}
\graphicspath{ {./images/} }
\usepackage{amssymb,epsfig,amsmath,amsthm,euscript,array} 
\usepackage[nosort]{cite}
\usepackage{pstricks}
\usepackage{polynom} 
\usepackage{color} 
\usepackage{tikz}
\usepackage{float}
\usepackage{caption}
\usepackage{subcaption}
\usepackage{caption}
\usepackage[T1]{fontenc}
\usepackage{arydshln} 
\usepackage{cancel} 
\usepackage{ulem} 
\usepackage{relsize} 
\usepackage{textcomp} 
\usepackage{url} 

\usepackage[unicode=true,pdfusetitle,
  bookmarks=true,bookmarksnumbered=true,bookmarksopen=true,
breaklinks=true,pdfborder={0 0 0},backref=false,colorlinks=true]{hyperref}

\hypersetup{ linkcolor=black, citecolor=black, urlcolor=black}

\makeatletter
\@addtoreset{equation}{section}
\makeatother

\newcounter{multieqs}

\newcommand{\be}{\begin{equation}}
\newcommand{\ee}{\end{equation}}

\def\lt{\tilde{\l}}

\newcommand{\bm}[1]{\mbox{\boldmath $#1$}}

\newcommand{\kslash}{k \!\!\! / }

\newcommand{\lslash}{l \!\! / }
\newcommand{\Pslash}{P \!\!\!\! / }

\newcommand{\islash}{i \!\!\! / }
\newcommand{\jslash}{j \!\!\! / }
\newcommand{\aslash}{a \!\!\! / }
\newcommand{\bslash}{{b \hspace{-6pt} \slash} }

\newcommand{\onslash}{1 \!\!\! / }
\newcommand{\twslash}{2 \!\!\!/ }
\newcommand{\thslash}{3 \!\!\!/ }
\newcommand{\foslash}{4 \!\!\! / }
\newcommand{\fislash}{5 \!\!\! / }

\newcommand{\mslash}{m \!\!\! / }

\def\bd{\begin{document}}
\def\ed{\end{document}}
\def\nn{\nonumber}
\def\bea{\begin{eqnarray}}
\def\eea{\end{eqnarray}}

\def\ab{(ijab)}
\def\ba{(ijba)}
\def\ijab{{\tr}_{+}(\islash\, \jslash\, \aslash \, \bslash)}
\def\ijba{{\tr}_{+}(\islash\, \jslash\, \bslash \, \aslash)}
\def\ijaP{{\tr}_{+}(\islash\, \jslash\, \aslash \, \Pslash)}
\def\ijPLa{{\tr}_{+}(\islash\, \jslash\, \Pslash_L \, \aslash)}
\def\ijaPL{{\tr}_{+}(\islash\, \jslash\, \aslash \, \Pslash_L)}
\def\ijPLza{{\tr}_{+}(\islash\, \jslash\, \Pslash_{L;z} \, \aslash)}
\def\ijaPLz{{\tr}_{+}(\islash\, \jslash\, \aslash \, \Pslash_{L;z})}
\def\ijPa{{\tr}_{+}(\islash\, \jslash\, \Pslash \, \aslash)}
\def\iaPb{{\tr}_{+}(\islash\, \aslash\, \Pslash \, \bslash)}
\def\ibPa{{\tr}_{+}(\islash\, \bslash\, \Pslash \, \aslash)}
\def\ijPmu{{\tr}_{+}(\islash\, \jslash\, \Pslash \, \mu)}
\def\ibmuP{{\tr}_{+}(\islash\, \bslash\, \mu \, \Pslash)}
\def\ibmua{{\tr}_{+}(\islash\, \bslash\, \mu \, \aslash)}
\def\iamub{{\tr}_{+}(\islash\, \aslash\, \mu \, \bslash)}
\def\jaPb{{\tr}_{+}(\jslash\, \aslash\, \Pslash \, \bslash)}
\def\ijmuP{{\tr}_{+}(\islash\, \jslash\, \mu \, \Pslash)}
\def\ijmum{{\tr}_{+}(\islash\, \jslash\, \mu \, \mslash)}
\def\ijmmu{{\tr}_{+}(\islash\, \jslash\, \mslash \, \mu)}
\def\ijmP{{\tr}_{+}(\islash\, \jslash\, \mslash \, \Pslash)}
\def\iabP{{\tr}_{+}(\islash\, \aslash\, \bslash \, \Pslash)}
\def\ijbP{{\tr}_{+}(\islash\, \jslash\, \bslash \, \Pslash)}
\def\jbPa{{\tr}_{+}(\jslash\, \bslash\, \Pslash \, \aslash)}
\def\ijPb{{\tr}_{+}(\islash\, \jslash\, \Pslash \, \bslash)}
\def\jbmua{{\tr}_{+}(\jslash\, \bslash\, \mu \, \aslash)}

\def\loablt{ {\tr}_{+}(\lslash_1\, \aslash \, \bslash\, \lslash_2)}

\def\ijlolt{{\tr}_{+}(\islash\, \jslash\, \lslash_1 \, \lslash_2)}
\def\ijltlo{{\tr}_{+}(\islash\, \jslash\, \lslash_2 \, \lslash_1)}
\def\ibloa{{\tr}_{+}(\islash\, \bslash\, \lslash_1 \, \aslash)}
\def\jaltb{{\tr}_{+}(\jslash\, \aslash\, \lslash_2 \, \bslash)}
\def\ialtb{{\tr}_{+}(\islash\, \aslash\, \lslash_2 \, \bslash)}
\def\bltloa{{\tr}_{+}(\bslash\, \lslash_2\, \lslash_1 \, \aslash)}
\def\jbloa{{\tr}_{+}(\jslash\, \bslash\, \lslash_1 \, \aslash)}
\def\ibPb{{\tr}_{+}(\islash\, \bslash\, \Pslash \, \bslash)}
\def\ijltb{{\tr}_{+}(\islash\, \jslash\, \lslash_2 \, \bslash)}

\def\ijloa{{\tr}_{+}(\islash\, \jslash\,  \lslash_1 \, \aslash)}
\def\ijblt{{\tr}_{+}(\islash\, \jslash\,  \bslash \, \lslash_2)}

\def\jakb{{\tr}_{+}(\jslash\, \aslash\, \kslash \, \bslash)}
\def\iakb{{\tr}_{+}(\islash\, \aslash\, \kslash \, \bslash)}

\def\tofo{{\tr}_{+}(\onslash\, \thslash\, \twslash \, \foslash)}
\def\foto{{\tr}_{+}(\onslash\, \thslash\, \foslash \, \twslash)}
\def\tofi{{\tr}_{+}(\onslash\, \thslash\, \twslash \, \fislash)}
\def\fito{{\tr}_{+}(\onslash\, \thslash\, \fislash \, \twslash)}

\def\lrangle#1#2{\langle #1\,#2\rangle}

\def\Li{{$\rm Li}_2$}
\def\eps{\epsilon}
\def\epsuv{{\epsilon_{\rm \mbox{\tiny UV}}}}
\let\bm=\bibitem
\let\la=\label

\def\npb#1#2#3{Nucl. Phys. {\bf{B#1}} #3 (#2)}
\def\plb#1#2#3{Phys. Lett. {\bf{#1B}} #3 (#2)}
\def\prl#1#2#3{Phys. Rev. Lett. {\bf{#1}} #3 (#2)}
\def\prd#1#2#3{Phys. Rev. {D \bf{#1}} #3 (#2)}
\def\cmp#1#2#3{Comm. Math. Phys. {\bf{#1}} #3 (#2)}
\def\cqg#1#2#3{Class. Quantum Grav. {\bf{#1}} #3 (#2)}
\def\nppsa#1#2#3{Nucl. Phys. B (Proc. Suppl.) {\bf{#1A}}#3 (#2)}
\def\ap#1#2#3{Ann. of Phys. {\bf{#1}} #3 (#2)}
\def\ijmp#1#2#3{Int. J. Mod. Phys. {\bf{A#1}} #3 (#2)}
\def\rmp#1#2#3{Rev. Mod. Phys. {\bf{#1}} #3 (#2)}
\def\mpla#1#2#3{Mod. Phys. Lett. {\bf A#1} #3 (#2)}
\def\jhep#1#2#3{J. High Energy Phys. {\bf #1} #3 (#2)}
\def\atmp#1#2#3{Adv. Theor. Math. Phys. {\bf #1} #3 (#2)}
\newcommand{\EQ}[1]{\begin{equation} #1 \end{equation}}
\newcommand{\AL}[1]{\begin{subequations}\begin{align} #1 \end{align}\end{subequations}}
\newcommand{\SP}[1]{\begin{equation}\begin{split} #1 \end{split}\end{equation}}
\newcommand{\ALAT}[2]{\begin{subequations}\begin{alignat}{#1} #2 \end{alignat}
  \end{subequations}}
\def\beqa{\begin{eqnarray}}
\def\eeqa{\end{eqnarray}}
\def\beq{\begin{equation}}
\def\eeq{\end{equation}}
\def\sst{\scriptscriptstyle}
\def\thetabar{\bar\theta}
\def\Tr{{\rm Tr}}
\def\one{\mbox{1 \kern-.59em {\rm l}}}
 \def\Nh{\hat{N}}
\newcommand{\half}{{\textstyle {1 \over 2}}}

\def\a{t}      \def\da{{\dott}}
\def\b{\beta}       \def\db{{\dot\beta}}
\def\c{\gamma}  \def\G{\Gamma}  \def\cdt{\dot\gamma}
\def\d{\delta}  \def\D{\Delta}  \def\ddt{\dot\delta}
\def\e{\epsilon}        \def\vare{\varepsilon}
\def\f{\phi}    \def\F{\Phi}    \def\vvf{\f}
\def\h{\eta}
\def\k{\kappa}
\def\l{\lambda} \def\L{\Lambda}
\def\m{\mu} \def\n{\nu}
\def\o{\omega}
\def\p{\pi} \def\P{\Pi}
\def\r{\rho}
\def\s{\sigma}  \def\S{\Sigma}
\def\t{\tau}
\def\th{\theta} \def\Th{\Theta} \def\vth{\vartheta}
\def\X{\Xeta}
\def\z{\zeta}
\def\de{\partial}
\def\cA{{\cal A}} \def\cB{{\cal B}} \def\cC{{\cal C}}
\def\cD{{\cal D}} \def\cE{{\cal E}} \def\cF{{\cal F}}
\def\cG{{\cal G}} \def\cH{{\cal H}} \def\cI{{\cal I}}
\def\cJ{{\cal J}} \def\cK{{\cal K}} \def\cL{{\cal L}}
\def\cM{{\cal M}} \def\cN{{\cal N}} \def\cO{{\cal O}}
\def\cP{{\cal P}} \def\cQ{{\cal Q}} \def\cR{{\cal R}}
\def\cS{{\cal S}} \def\cT{{\cal T}} \def\cU{{\cal U}}
\def\cV{{\cal V}} \def\cW{{\cal W}} \def\cX{{\cal X}}
\def\cY{{\cal Y}} \def\cZ{{\cal Z}}
\def\ua{\underline{t}}
\def\ub{\underline{\phantom{t}}\!\!\!\beta}
\def\uc{\underline{\phantom{t}}\!\!\!\gamma}
\def\um{\underline{\mu}}
\def\ud{\underline\delta}
\def\ue{\underline\epsilon}
\def\una{\underline a}\def\unA{\underline A}
\def\unb{\underline b}\def\unB{\underline B}
\def\unc{\underline c}\def\unC{\underline C}
\def\und{\underline d}\def\unD{\underline D}
\def\une{\underline e}\def\unE{\underline E}
\def\unf{\underline{\phantom{e}}\!\!\!\! f}\def\unF{\underline F}
\def\unm{\underline m}\def\unM{\underline M}
\def\unn{\underline n}\def\unN{\underline N}
\def\unp{\underline{\phantom{a}}\!\!\! p}\def\unP{\underline P}
\def\unq{\underline{\phantom{a}}\!\!\! q}
\def\unQ{\underline{\phantom{A}}\!\!\!\! Q}
\def\unH{\underline{H}}
\def\As {{A \hspace{-6.4pt} \slash}\;}
\def\bs {{b \hspace{-6.4pt} \slash}\;}
\def\Ds {{D \hspace{-6.4pt} \slash}\;}
\def\ds {{\del \hspace{-6.4pt} \slash}\;}
\def\ss {{\s \hspace{-6.4pt} \slash}\;}
\def\ks {{ k \hspace{-6.4pt} \slash}\;}
\def\ps {{p \hspace{-6.4pt} \slash}\;}
\def\pas {{{p_1} \hspace{-6.4pt} \slash}\;}
\def\pbs {{{p_2} \hspace{-6.4pt} \slash}\;}
\def\Ps {{P \hspace{-6.4pt} \slash}\;}
\def\Qs {{Q \hspace{-6.4pt} \slash}\;}
\def\Fh{\hat{F}}
\def\Vh{\hat{V}}
\def\Xh{\hat{X}}
\def\ah{\hat{a}}
\def\xh{\hat{x}}
\def\yh{\hat{y}}
\def\ph{\hat{p}}
\def\xih{\hat{\xi}}
\def\psit{\tilde{\psi}}
\def\Psit{\tilde{\Psi}}
\def\tht{\tilde{\th}}
\def\lt{\tilde{\lambda}}
\def\hl{\hat{\lambda}}
\def\hlt{\hat{\tilde{\lambda}}}
\def\llt{\tilde{l}}
\def\At{\tilde{A}}
\def\Qt{\tilde{Q}}
\def\Rt{\tilde{R}}
\def\Nt{\tilde{N}}
\def\at{\tilde{a}}
\def\st{\tilde{s}}
\def\ft{\tilde{f}}
\def\pt{\tilde{p}}
\def\qt{\tilde{q}}
\def\vt{\tilde{v}}
\def\nt{\tilde{n}}
\def\delb{\bar{\partial}}
\def\bz{\bar{z}}
\def\bD{\bar{D}}
\def\bB{\bar{B}}
\def\bk{{\bf k}}
\def\bl{{\bf l}}
\def\bp{{\bf p}}
\def\bq{{\bf q}}
\def\br{{\bf r}}
\def\bx{{\bf x}}
\def\by{{\bf y}}
\def\bR{{\bf R}}
\def\bV{{\bf V}}
\def\d{\delta}\def\D{\Delta}\def\ddt{\dot\delta}
\def\pa{\partial} \def\del{\partial}
\def\xx{\times}
\def\uno{\mbox{1 \kern-.59em {\rm l}}}
\def\trp{^{\top}}
\def\inv{^{-1}}
\def\dag{{^{\dagger}}}
\def\pr{^{\prime}}
\def\lan{\langle}
\def\ran{\rangle}
\def\rar{\rightarrow}
\def\lar{\leftarrow}
\def\lrar{\leftrightarrow}
\newcommand{\0}{\,\!}      
\def\one{1\!\!1\,\,}
\def\im{\imath}
\def\jm{\jmath}
\newcommand{\tr}{\mbox{tr}}
\newcommand{\slsh}[1]{/ \!\!\!\! #1}
\def\vac{|0\rangle}
\def\lvac{\langle 0|}
\def\hlf{\frac{1}{2}}
\def\ove#1{\frac{1}{#1}}
\def\Box{\square}
\def\ZZ{\mathbb{Z}}
\def\CC#1{({\bf #1})}
\def\bcomment#1{}
\def\bfhat#1{{\bf \hat{#1}}}
\def\VEV#1{\left\langle #1\right\rangle}
\newcommand{\ex}[1]{{\rm e}^{#1}} \def\ii{{\rm i}}
\def\rr{{\rm r}} \def\rs{{\rm s}}\def\rv{{\rm v}}
\def\ri{{\rm i}}\def\rj{{\rm j}}
\newcommand{\lrbrk}[1]{\left(#1\right)}
\newcommand{\sfrac}[2]{{\textstyle\frac{#1}{#2}}}
\def\Li{{\rm Li}_2}
\DeclareMathOperator{\dif}{d \!}

\newcommand{\ang}[1]
{
  \langle #1 \rangle
}

\newcommand{\sqr}[1]
{
  [ #1 ]
}
\newcommand{\tlambda}
{
  \tilde{\lambda}
}

\newcommand{\diffd}{\mathrm{d}}

\newcommand{\mmbox}[3]{\makebox[#1][#2]{#3}}

\newcommand*\circled[1]{\tikz[baseline=(char.base)]{
            \node[shape=circle,draw,inner sep=0.5pt] (char) {\tiny #1};}}

\font\mybb=msbm10 at 12pt
\def\bb#1{\hbox{\mybb#1}}

\font\myBB=msbm10 at 18pt
\def\BB#1{\hbox{\myBB#1}}

\newcommand\note[1]{\mbox{}\marginpar{\footnotesize\raggedright\hspace{0pt}\color{blue}\emph{#1}}}

\newtheorem{theorem}{Theorem}
\newtheorem{lemma}{Lemma}
\newtheorem{proposition}{Proposition}

\theoremstyle{definition}
\newtheorem{corollary}{Corollary}
\newtheorem{definition}{Definition}
\newtheorem{conjecture}{Conjecture}
\newtheorem{hypothesis}{Hypothesis}
\newtheorem{remark}{Proof}
\newtheorem{rmk}{Remark}
\newtheorem{example}{Example}

\begin{document}
\begin{flushright}
\end{flushright}

\begin{center}

{\Large \bf   A Generalisation of Ramanujan's (back of the envelope) Method for Divergent Series}

\vspace{12pt}

{\mbox {\bf Patrick~J.~Burchell}}%

\end{center}

\begin{abstract}
Ramanujan derived the well known divergent-sum of integers in more than one way. We generalise the informal method to higher powers of the Riemann zeta function through a study of the Eulerian numbers in particular. Within the context of additive combinatorics a heuristic approach that unifies generating series and difference matrices is presented. 
\end{abstract}

\setcounter{page}{1}
\thispagestyle{empty}

\section{Introduction}
Within this paper we evaluate a seemingly informal method, the outcome of which was presented to Hardy by Ramanujan in his letter of 16 January 1913, and which Berndt details in typeset form within his presentation of Ramanujan's Notebook 2.\!\cite{ramanujan_watson_1}\!\cite{ramanujan_berndt} In fact it appears verbatim in Notebook 1 and Notebook 2.%
\footnote{The sketch here presented is from Notebook 1.}\!\cite{ramanujan_1}\!\cite{ramanujan_2} Intriguingly Berndt observes that it is peculiar that Ramanujan utilises this method when he has his \textit{own} abbreviated Euler--Maclaurin sum formula immediately to hand.\!\cite{ramanujan_berndt} Thus within this paper a generalisation of the steps is presented in order to achieve a broader context within generating series. \\
\begin{figure}[h!]
\centering\includegraphics[scale=0.91]{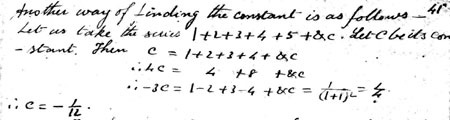}
\end{figure}
\\
Critically a general relation, in terms of a closed-form ordinary generating function, at Equation~\ref{eq:RamGen}, is presented which provides the impetus to this paper. The standard $\zeta(s)$ to $\eta(s)$ comparison and transition features at Proposition~\ref{prp:Upsilon} but whilst the latter is known it is the author's understanding that the former is not. Indeed a particular theme of the closed form algorithm is the appearance of the Eulerian numbers. Lemma~\ref{lem:BernoulliEulerian} exploits a fascinating relation between them and the Bernoulli numbers to arrive at the well known relation $\zeta(-s)=-\frac{B_{s+1}}{s+1}$ for all strictly positive values of $s$. Finally Lemma~\ref{lem:GeneratingSeriesDifference} unifies difference matrices and generating series in an approach that distinctly favours the backward difference operator.
\section{Ramanujan's Method \& the Ordinary Generating Series} 
We begin by representing Ramanujan's method in terms of generating series. The penultimate step of Ramanujan's note indeed suggests this in no uncertain terms. Where for instance he writes $1/(1+1)^2$. This is the classic step of substituting unity into the variable within a formal power series. Hence we have
\beq
  1-2z^2+3z^3-...=\frac{1}{(1+z)^2}
\, ,
\eeq
with the closed-form on the right, and the open-form left. As per the second to last line of Ramanujan's note.
\begin{example}
\label{ex:Power1}
If we select Ramanujan's example $f(x)=(x+1)$, in functional terms, we can reproduce it in terms of a generating series\cite{wilf} as
{
\[ 
  \begin{tabular}{ccccccc}
       &  $1z^0$ & $+2z^1$   & $+3z^2$      & $+4z^3$  & $+5z^4$ & ...  \\  
   $-($ &                  & $\phantom{+}4z^1$ &                        & $+8z^3$ & & \phantom{  } ... $)$  \\
   \hline
       & $1z^0$ & $-2z^1$   & $+3z^2$ & $-4z^3$ & $+5z^4$ & ...
  \end{tabular}
\]
}
so that in terms of the closed-form generating series we have
\beq
\label{eq:avoid_division_zero}
   \frac{1}{(1-z)^{2}}-\frac{2^{1+1}z}{(1-z)^{2}(1+z)^{2}}=\frac{1}{(1+z)^{2}}
\, .
\eeq
\end{example}
\begin{rmk}
It's important to note that we have cued off three things: (a) $f(x)=x+1=\{1,2,3,...\}$, (b) the final key step $1/(1+1)^2$ which matches $1/(1+z)^2$ not $z/(1+z)^2$, (c) the open-form of the generating series, top and bottom. So that we start with the powers, and sum to the alternating powers.
\end{rmk}
\begin{rmk}
\label{rmk:avoid_division_zero}
If we set aside deeper considerations of the \textit{disk of convergence} it is clear that under the ordinary generating series, more formally defined under Theorem~\ref{thm:MAIN}, there is one particular instance where substitution of $z=1$ does not work at all. In fact it is the very generating series that produces the polynomials given as $g(z)=\frac{t_0+t_1z+...+t_sz^s}{(1-z)^{s+1}}$. This is a challenge that must be overcome within the example. Indeed, on the left hand side of Equation~\ref{eq:avoid_division_zero} the two fractions combine to achieve just that. This same obstacle is encountered as one extrapolates beyond the first power. It is not unfeasible to argue that overcoming it might have caused Ramanujan to value this algorithm beyond the obvious result. 
\end{rmk}
\begin{rmk}
Obscured (if not quite hidden) in plane sight is the technique to extend this method whereby we observe that $\frac{1}{(1-z)^{2}}-\big[\frac{1}{(1-z)^{2}}-\frac{1}{(1+z)^{2}}\big]=\frac{1}{(1+z)^{2}}$ is a restatement of Equation~\ref{eq:avoid_division_zero}.
\end{rmk}
Clearly the target, the alternating values, must be a result that one knows, and the middle sequence must be codifiable in terms of Ramanujan's $C$. None of which can be assumed. To give the reader some confidence however we will step up to the 3rd powers.
\begin{example}
\label{ex:Power3}
For the 3rd powers we have
\[ 
  \begin{tabular}{ccccccc}
       &  $1z^0$ & $+8z^1$   & $+27z^2$      & $+64z^3$  & $+125z^4$ & ...  \\  
   $-($ &                  & \phantom{$.$}$16z^1$ &                        &\phantom{$.$}$128z^3$ &    & \phantom{  } ... $)$ \\
   \hline
       & $1z^0$ & $-8z^1$   & $+27z^2$ & \phantom{.}$-64z^3$ & $+125z^4$ & ...
  \end{tabular}
\]
which equates to
\beq
   \frac{1+4z+z^2}{(1-z)^{3+1}}-\frac{2^{3+1}(z+4z^3+z^5)}{(1-z)^{3+1}(1+z)^{3+1}}=\frac{1-4z+z^2}{(1+z)^{3+1}}
\, ,
\eeq
in closed-form. The pattern in the numerator term(s) does indeed repeat and is given by the Eulerian numbers. Naturally Ramanujan's manipulation of the constant, designated by $C$, goes from $C-2^{1+1}C$ to $C-2^{3+1}C$. And quite remarkably the math works out. Here we have
\beq
  C-2^4C=\frac{1-4+1}{(1+1)^4}
\, ,
\eeq
where $C=1/120=\zeta(-3)$. 
\end{example}
\section{A General Relation}
\begin{theorem}
\label{thm:MAIN}
Let $\mathbb{R}[[z]]$ define the ring of formal power series in the variable $z$ over $\mathbb{R}$.%
\footnote{This proof will not lean upon subtle topological properties of formal power series. Which is not to say that it will not highlight them.} For the formal power series that represents the powers in $x,s \in \mathbb{Z}$ given as $f(x)=(x+1)^s$, $x \geq 0$, $s \geq 1$, it is held that
\beq
\label{eq:RamGen}
  \frac{\Upsilon_{\!s}}{(1-z)^{s+1}}-2^{s+1} \! \cdot \! \frac{\Lambda_s}{(1-z)^{s+1}(1+z)^{s+1}}=\frac{\xi_{s}}{(1+z)^{s+1}}
\, .
\eeq
Where
\beq
\label{eq:OffsetEulerian}
  \Upsilon_{\!s}=\sum_{k=1}^{s} 
\sbox0{$\begin{array}{c}s\\k\end{array}$}
\resizebox{.4\width}{.8\ht0}{$\Bigg\langle$}
\!
\usebox{0}
\!
\resizebox{.4\width}{.8\ht0}{$\Bigg\rangle$}\,z^{k-1}
\, ,
\eeq
and
\beq
\label{eq:Lambda}
  \Lambda_s= \sum_{k=1}^{s} 
\sbox0{$\begin{array}{c}s\\k\end{array}$}
\resizebox{.4\width}{.8\ht0}{$\Bigg\langle$}
\!
\usebox{0}
\!
\resizebox{.4\width}{.8\ht0}{$\Bigg\rangle$}\,z^{2k-1}
\, ,
\eeq
and
\beq
\label{AlternatingPowersEta}
  \xi_{s}=\sum_{k=1}^{s} 
\sbox0{$\begin{array}{c}s\\k\end{array}$}
\resizebox{.4\width}{.8\ht0}{$\Bigg\langle$}
\!
\usebox{0}
\!
\resizebox{.4\width}{.8\ht0}{$\Bigg\rangle$}\,z^{k-1}(-1)^{k-1}
\, 
\eeq
respectively. Further, it is held that this relation generalises Ramanujan's method where we set firstly $\Upsilon_{\!s}\! \cdot \!(1-z)^{-s-1}=g(z)$ and then $\Lambda_s\! \cdot \!(1-z)^{-s-1}(1+z)^{-s-1}=h(z)$ to produce
\beq
\label{eq:EulerBernoulli1}
  g(z)-2^{s+1}h(z)=\frac{\xi_{s}}{(1+z)^{s+1}}
.
\eeq
We then substitute $z=1$, directly in the right hand side, but indirectly, through the open-form of $g(z)$ and $h(z)$ to give $C-2^{s+1}C$ and
\beq
\label{eq:EulerBernoulli2}
  C=\frac{\xi_s}{2^{s+1}(1-2^{s+1})}
\, .
\eeq
To produce $C=\zeta(-s) \,\, \forall \,\, s \geq 1$.%
\footnote{This covers powers equal to and greater than the integers (at unity) within Ramanujan's sketch.}
\end{theorem}
\begin{rmk}
Observe that the Eulerian numbers\cite{euler} appear in the closed-form of the generating series that captures the powers. If we wish to represent $f(x)=(x+1)^3$, $x \geq 0$, in terms of its open-form generating series, we have $g(z)=1+8z+27z^2+...$ whilst the closed-form is given as
\beq
  g(z)=\frac{1+4z+z^2}{(1-z)^4}
\, .
\eeq
Thus, for the purpose of this paper the Eulerian number triangle is converted into its generating series form (sans the denominator) by Equation~\ref{eq:OffsetEulerian}. Given a sample of the Eulerian number triangle\cite{knuth}
{
\[
\begin{tabular}{>{\centering\arraybackslash}p{4mm}>{\centering\arraybackslash}p{4mm}>{\centering\arraybackslash}p{4mm}>{\centering\arraybackslash}p{4mm}>{\centering\arraybackslash}p{4mm}>{\centering\arraybackslash}p{4mm}>{\centering\arraybackslash}p{4mm}}
&&&&&&\\
1&&&&&&\\
1&1&&&&&\\
1&4&1&&& &\\
1&11&11&1&&&\\
1&26&66&26&1&&\\
..&..&..&..&..&..&\\
\end{tabular}
\label{tab:eulerian}
\]
}
we can read off the cubes as
\beq
\label{eq:OffsetEulerian3}
  \Upsilon_{\!3}=\sum_{k=1}^{3} 
\sbox0{$\begin{array}{c}3\\k\end{array}$}
\resizebox{.4\width}{.8\ht0}{$\Bigg\langle$}
\!
\usebox{0}
\!
\resizebox{.4\width}{.8\ht0}{$\Bigg\rangle$}\,z^{k-1}
\, ,
\eeq 
which gives us $1+4z+z^2$. It follows that within the large brackets $s$ captures the power and row, which is tallied from $1$, and $k$ similarly represents the column tallied from unity.
\end{rmk}
Whilst we assume the powers, and hence the $\Upsilon_{\!s}$ term,\cite{knuth} as a given, the more exotic $\xi_s$ and $\Lambda_s$ terms will be clearly demonstrated below. They will however draw numbers from the Eulerian triangle in a similar manner.
\section{The Alternating Powers}
\begin{proposition}
The alternating powers $f(x)=(x+1)^s(-1)^x$, $s, x \in \mathbb{Z}$, $x \geq 0$, $s \geq 1$ are captured in generating series closed-form by $k(z)=\xi_{s}\! \cdot \!(1+z)^{-s-1}$, where $\xi_s$ is defined by Equation~\ref{AlternatingPowersEta}.
\end{proposition}
\begin{proof}
Now we know that the numerator terms of the closed-form generating series for powers are given by the Eulerian numbers. Specifically, this means that $\Upsilon_{\!s}\! \cdot \!(1-z)^{-s-1}=1^s+2^sz+3^sz^2+4^sz^3+...$ with $\Upsilon_{\!s}$ defined at Equation~\ref{eq:OffsetEulerian}. Now, in the form of Equation~\ref{eq:OffsetEulerian}, we need only substitute $-z$ to obtain $\xi_{s}\! \cdot \!(1+z)^{-s-1}=1^s-2^sz+3^sz^2-4^sz^3+...$ where $\xi_s$ is given under Equation~\ref{AlternatingPowersEta}.
\end{proof}
\begin{proposition}
\label{prp:Upsilon}
It is held that 
\beq
  \Upsilon_{\!s}\! \cdot \!(1-z)^{-s-1}-\xi_s\! \cdot \!(1+z)^{-s-1}=h(z)\! \cdot \!2^{s+1}
\, ,
\eeq
where $h(z)$ is given in open-form as 
\beq
 h(z)=1^sz+2^sz^3+3^sz^5+4^sz^7+5^sz^9+...
\, 
\eeq
with closed form to be determined.
\end{proposition}
\begin{proof}
We assess $\Upsilon_{\!s}\! \cdot \!(1-z)^{-s-1}-\xi_s\! \cdot \!(1+z)^{-s-1}$ in its open-form.
{
\[ 
  \begin{tabular}{cccccccc}
       &  $1^{s}$ & $+2^{s}z$   & $+3^{s}z^2$      & $+4^{s}z^3$  & $+5^{s}z^4$ & $+6^{s}z^5$      &  \phantom{$)$} \,\, ...  \\  
   $- \,\, ($     &  $1^{s}$ & $-2^{s}z$   & $+3^{s}z^2$      & $-4^{s}z^3$  & $+5^{s}z^4$ & $-6^{s}z^5$       & $)$ \,\, ...  \\  
   \hline
       &          & $2\! \cdot \!2^{s}z$   &   $+$\phantom{$3^{s}$}    & $2\! \cdot \!4^{s}z^3$ &   $+$\phantom{$5^{s}$}       &  $2\! \cdot \!6^{s}z^5$          &  \phantom{$)$} \,\, ...
  \end{tabular}
\]
}
It follows by inspection that $\Upsilon_{\!s}\! \cdot \!(1-z)^{-s-1}-\xi_s\! \cdot \!(1+z)^{-s-1}=h(z)\! \cdot \!2^{s+1}$.%
\footnote{A more familiar formula tells us that $\eta(s)=(1-s^{1-s})\zeta(s)$. If we substitute $-s$ and rearrange we obtain $\zeta(-s)=\eta(-s)(1-2^{s+1})^{-1}$. See Example~\ref{ex:Power3}.}
\end{proof}
\begin{rmk}
This step tells us (a) that $h(z)$ exists, and (b) that it replicates the governing series $g(z)=\Upsilon_{\!s}\! \cdot \!(1-z)^{-s-1}$. This is the case \textit{specifically} when we substitute $z=1$ within the open-form, so that we obtain $g(z)=h(z)$.%
\footnote{Naturally this is the contentious step of Ramanujan's sketch.} And these facts hold independently of whether the closed form of $h(z)$ is $h(z)=\Lambda_s\! \cdot \!(1-z)^{-s-1}(1+z)^{-s-1}$, as we have postulated within Equation~\ref{eq:RamGen}, and Equation~\ref{eq:Lambda}.
\end{rmk}
\section{The Bernoulli and Eulerian Numbers}
\begin{lemma}
\label{lem:BernoulliEulerian}
The step from Equation~\ref{eq:EulerBernoulli2} to the relation $C=\zeta(-s) \,\, \forall \,\, s \geq 1$, at Theorem~\ref{thm:MAIN}, can be deduced from a known relation between the Eulerian and Bernoulli numbers,\cite{luschny} given as
\beq
\label{eq:EulerBernoulliKnown}
  \sum_{k=1}^{s} 
\sbox0{$\begin{array}{c}s\\k\end{array}$}
\resizebox{.4\width}{.8\ht0}{$\Bigg\langle$}
\!
\usebox{0}
\!
\resizebox{.4\width}{.8\ht0}{$\Bigg\rangle$}\,(-1)^{k-1}=-\beta_{s+1}\! \cdot \!(2^{s+1}-4^{s+1})
\, ,
\eeq
with $\beta_{s+1}=\frac{B_{s+1}}{s+1}$, the divided Bernoulli numbers.
\end{lemma}

\begin{proof}
If one substitutes $z=1$ under Equation~\ref{AlternatingPowersEta} then $\xi_s$ reduces to the alternating sum of the Eulerian numbers 
\beq
\label{eq:ReducedXi}
  \xi_s=\sum_{k=1}^{s} 
\sbox0{$\begin{array}{c}s\\k\end{array}$}
\resizebox{.4\width}{.8\ht0}{$\Bigg\langle$}
\!
\usebox{0}
\!
\resizebox{.4\width}{.8\ht0}{$\Bigg\rangle$}(-1)^{k-1}
\, .
\eeq
Thus combining with Equation~\ref{eq:EulerBernoulliKnown} we have $\xi_s=-\beta_{s+1}\! \cdot \!(1-2^{s+1})2^{s+1}$. Finally, if we sub into Equation~\ref{eq:EulerBernoulli2} we find that
\beq
  C=-\beta_{s+1}
\, ,
\eeq
which is a standard formula for the Riemann Zeta function: $\zeta(-s)=-\beta_{s+1}$ for all $s \geq 0$, $s \in \mathbb{Z}$.\cite{riemann_edwards}
\end{proof}
\section{Deconvolution}
\begin{lemma}
\label{lem:Upsilon}
At Proposition~\ref{prp:Upsilon} we proved $\Upsilon_{\!s}\! \cdot \!(1-z)^{-s-1}-\xi_s\! \cdot \!(1+z)^{-s-1}=h(z)\! \cdot \!2^{s+1}$ with the open-form of $h(z)$ given by $h(z)=1^sz+2^sz^3+3^sz^5+4^sz^7+5^sz^9+...\,\,$. We now hold that the closed form is such that $h(z)=\Lambda_s\! \cdot \!(1-z)^{-s-1}(1+z)^{-s-1}$ where $\Lambda_s$ is given at Equation~\ref{eq:Lambda}.
\end{lemma}
The proof of this Lemma, which can be found at Appendix~\ref{sec:OGS}, completes the proof of Theorem~\ref{thm:MAIN}.
\begin{rmk}
The terms deconvolution and convolution are naturally context specific. Here we utilise the word deconvolution (in its first instance) to mean the reversal of the difference operation. Although it has an easy analogue as the partial sum this falls away when one looks at denominator terms of greater complexity, whereupon the general deconvolution matrix presented within Appendix~\ref{sec:OGS} becomes a more practicable alternative. 
\end{rmk}
\section{Final Worked Example}
Note that for even values of $s$ it is easy to see that $C$ reduces to zero since alternating sums of the Eulerian numbers at even rows cancel.
\begin{example}
\label{ex:Power5}
If we select as example $f(n)=(n+1)^5$ we can reproduce the left hand side of \ref{eq:RamGen} as
\beq
\label{eq:Power5}
   \frac{1+26z+66z^2+26z^3+z^4}{(1-z)^{s+1}}-\frac{64(z+26z^3+66z^5+26z^7+z^9)}{(1-z)^{s+1}(1+z)^{s+1}}
\, ,
\eeq
which reduces to $(1-26z+66z^2-26z^3+z^4)(1+z)^{-s-1}$ per Equation~\ref{eq:RamGen}. Next we resolve as per the instructions so that we have
\beq
  C-2^6\! \cdot \!C=\frac{1-26+66-26+1}{(1+1)^6}
\, .
\eeq
Which reduces to $C=-1/252=\zeta(-5)$ as required. What's more we have the open-form sum
{
\[ 
  \begin{tabular}{ccccccc}
       &  $1z^0$ & $+32z^1$   & $+243z^2$      & $+1024z^3$  & $+3125z^4$ & ...  \\  
   $-($ &                  & \phantom{$+$}$64z^1$ &                        & $+2048z^3$ &                 & \phantom{...}... $)$\\
   \hline
       & $1z^0$ & $-32z^1$   & $+243z^2$ & $-1024z^3$ & $+3125z^4$ & ...
  \end{tabular}
\]
}
and $h(z)$, within the second term of Equation~\ref{eq:Power5}, given by
\beq
  h(z)=z+32z^3+243z^5+1024z^7+3125z^9+...
\,\,\,\, .
\eeq
\end{example}
\section{Conclusion}
It has been shown that the method given within Theorem~\ref{thm:MAIN} will produce $\zeta(-s)$, $s \in \mathbb{Z}^+$ in all instances, and that it mirrors and in fact generalises Ramanujan's technique. Within Theorem~\ref{thm:MAIN}, and following Ramanujan's example, $C$ is equated to two seemingly distinct objects: $g(z)$ and $h(z)$. However, one can see that the open-form of both generating series does not differ when one substitutes $z=1$. Although this does not address deeper topological concerns over the validity of that substitution it does appraise the consequences of that substitution within the context of a general method.
\section{Acknowledgments}
The author is very grateful to Bruce Berndt Professor Emeritus at the University of Illinois and to Damian R\"ossler of the University of Oxford for valuable insights and feedback.
\appendix
\section{Deconvolution of Ordinary Generating Series}
\label{sec:OGS}
This section gives a proof of Lemma~\ref{lem:Upsilon} in a manner that (it is felt) provides a particular heuristic value. Fundamentally one wishes to show that
\beq
  1^sz+2^sz^3+3^sz^5+4^sz^7+5^sz^9+...=\frac{\Lambda_s}{(1-z)^{s+1}(1+z)^{s+1}}
\, 
\eeq
with $\Lambda_s$ given under Equation~\ref{eq:Lambda}.
\begin{example} 
\label{ex:EULERIAN}
Thus we eschew the habit of placing difference terms between parent values, move our focus away from the common difference term, and towards the first finite difference sequence. The left-most column is $j=0$ and $t_{0,j}$ is the first sequence: the integers, and the squares respectively.
\[ 
  \begin{tabular}{ccccccc}
  0 & 1 & 2 & 3 & 4 & 5 & ...\\
  0 & 1 & 1 & 1 & 1 & 1 & ...\\
  0 & 1 & 0 & 0 & 0 & 0 & ...\\
   &  &  &  &  &  &  \\
  \end{tabular}
\,\,\,\,\,\,\,
  \begin{tabular}{ccccccc}
  0 & 1 & 4 & 9 & 16 & 25 & ... \\
  0 & 1 & 3 & 5 & 7 & 9 & ... \\
  0 & 1 & 2 & 2 & 2 & 2 & ... \\
  0 & 1 & 1 & 0 & 0 & 0 & ... \\
  \end{tabular}
\]
And so too for the cubes and the powers of four.
\[ 
\, \, \, \,
  \begin{tabular}{ccccccc}
  0 & 1 & 8 & 27 & 64 & 125 & ...\\
  0 & 1 & 7 & 19 & 37 & 61 & ...\\
  0 & 1 & 6 & 12 & 18 & 24 & ...\\
  0 & 1 & 5 & 6 & 6 & 6 & ...\\
  0 & 1 & 4 & 1 & 0 & 0 & ...\\
   &  &  &  &  &  &  \\
  \end{tabular}
\,\,\,\,\,\,\,
  \begin{tabular}{ccccccc}
  0 & 1 & 16 & 81 & 256 & 625 & ... \\
  0 & 1 & 15 & 65 & 175 & 369 & ... \\
  0 & 1 & 14 & 50 & 110 & 194 & ... \\
  0 & 1 & 13 & 36 & 60 & 84 & ... \\
  0 & 1 & 12 & 23 & 24 & 24 & ... \\
  0 & 1 & 11 & 11 & 1 & 0 & ... \\
  \end{tabular}
\]
Quite remarkably we have discovered the Eulerian numbers within the final difference sequence.
\end{example}
\begin{rmk}
\label{rmk:difference_op}
It is likely that this goes unheralded (in so far as it does) because one is, as a rule, more concerned with the repeating term or the line of the binomial transform. Which is to say that we would not have discovered these terms if we had performed a forward difference operation. Note that one can start from the bottom row and perform a partial sum to obtain all the sequences above, one a time. This is a special type of deconvolution. See Proposition~\ref{prop:PartialSum}.
\end{rmk}
\subsection{The Working Matrix}
\begin{definition}
\label{def:WorkingMatrix}
Define a working matrix with $l$ rows, counting from zero, as
\begin{equation}
t_{ij}
=
\left(
\begin{array}{rrrrrrr}
\dotsb & 0     &0 & t_{0,0}      & t_{0,1}    &t_{0,2}    & \dotsb\\
\dotsb & 0     & 0& t_{1,0}   &t_{1,1}      &t_{1,2}     &  \dotsb\\
\dotsb & 0     & 0&t_{2,0}   & t_{2,1}     &t_{2,2}     &  \dotsb\\
\reflectbox{$\ddots$} & \phantom{-}\vdots & \phantom{-}\vdots & \phantom{-}\vdots & \phantom{-}\vdots & \phantom{-}\vdots & \ddots\\
\dotsb & 0 &0 & t_{l,0} &t_{l,1} & t_{l,2} & \dotsb \\
\end{array}
\right) 
\, .
\end{equation}
The nature of the terms $t_k$ will be elucidated within the discussion to follow.
\end{definition}
Now combine Example~\ref{ex:EULERIAN} and Definition~\ref{def:WorkingMatrix}.
\begin{definition}
\label{def:DifferenceSequenceMatrix}
Define a non-finite sequence $(t_n)_{n \in \mathbb{N}}$, with $t_n \in \mathbb{R}$. Then define a difference sequence matrix $t_{i,j}$ that takes as the first input row $t_{0,n}=t_n$. Allow a suitable number of terms with $j\leq0$, all zero, and apply the formula $t_{i,j}=t_{i-1,j}-t_{i-1,j-1}$. When one accumulates these differences across multiple rows, say at the $l^{th}$ row, one obtains the formula
\beq
  t_{l,j}=\sum_{k=0}^{j} (t_{0,j-k})(-1^{k})\binom{l}{k}
\, .
\eeq
The matrix, for current purposes, ends at the first finite difference sequence. Naturally, if there is none, then the input sequence is not that of a univariate polynomial.
\end{definition}
\begin{rmk}
Note that, by definition, when a convolution (think difference calculation) extends beyond the first term, it will always return zero. Informally we refer to this as the \textit{off piste returns zero} rule. Holding this in mind allows us to significantly simplify formulae without loss of rigour.
\end{rmk}
\begin{definition}
\label{def:DSK}
Given a difference sequence matrix per Definition~\ref{def:DifferenceSequenceMatrix}, with $l$ rows, we formalise the first finite difference sequence as the \textit{difference sequence key}. Utilising the Eulerian number's angled braces for continuity we have $\langle t_{l,0},t_{l,1},..., t_{l,w} \rangle$ or preferably $\langle t_0,t_1,..., t_w \rangle_l$. 
\end{definition}
\begin{example}
For the cubes (including zero) their \textit{difference sequence key} is given by $\langle 0,1,4,1 \rangle_{3+1}$.
\end{example}
\begin{rmk}
The d.s.k. isolates and exploits particular values that spill out of a difference matrix but are lost under a forward difference analysis.
\end{rmk}
\begin{proposition}
\label{prop:PartialSum}
The partial sum of the $l$th row of a difference matrix $t_{i,j}$ given per Definition~\ref{def:DifferenceSequenceMatrix} exactly reverses the difference operation, where the partial sum is given as
\beq
\label{eq:PartialSum}
  t_{l-1,j}=\sum_{k=0}^j t_{l,k}
\, .
\eeq 
\end{proposition}
\begin{proof}
Hence from Definition~\ref{def:DifferenceSequenceMatrix} we build a difference matrix with the final two rows given below.
\[ 
  \begin{tabular}{cccccc}  
  $t_{l-1,0}$ & $t_{l-1,1}$ & $t_{l-1,2}$ & $t_{l-1,3}$ & $t_{l-1,4}$ & ...\\
  $t_{l,0}$ & $t_{l,1}$ & $t_{l,2}$ & $t_{l,3}$ & $t_{l,4}$ & ...\\
  \end{tabular}
\]
And which, for aesthetic reasons, we replace with variables.
\[ 
  \begin{tabular}{rrrrrr}  
    A & B & C & D & E & ...\\
    a & b & c & d & e & ...\\
  \end{tabular}
\]
Now, from the defined difference operation we have: $A=a$, $B-A=b$, $C-B=c$, and $D-C=d$ etc. This can be restated as $A=a$, $B=b+a$, $C=c+b+a$, and $D=d+c+b+a$ etc., thus producing the partial sum algorithm at \ref{eq:PartialSum}.
\end{proof}
What this means is that the matrices can be constructed as difference operations on the terms of the input sequence, or as partial sums operating on the terms of the difference key.
\begin{definition}
Given a partial sum matrix per Definition~\ref{def:DifferenceSequenceMatrix} and Proposition~\ref{prop:PartialSum}, and tallying the d.s.k. row from zero, then the sequence found at row $w+1$ is called the midline sequence.%
\footnote{Within Example~\ref{ex:EULERIAN} these are the input rows.}
\end{definition}
\subsection{The Ordinary Generating Series \& Difference Operations}
\begin{lemma}
\label{lem:GeneratingSeriesDifference}
Let $\mathbb{R}[[z]]$ define the ring of formal power series in the variable $z$ over $\mathbb{R}$. Given a difference matrix per Definition~\ref{def:DifferenceSequenceMatrix} with $(t_n)_{n \in \mathbb{N}}$, $t_n \in \mathbb{R}$, and where $t_{0,n}=t_n$, we define a complementary ordinary generating series such that
\beq
  g_0(z)=t_{0,0}+t_{0,1}z+t_{0,2}z^2+...
\,\,\, .
\eeq
It is held that the operation $g_0(z)*(1-z)$ is term-wise equivalent to the difference operation under Definition~\ref{def:DifferenceSequenceMatrix}. Thus, for a matrix with $l$ rows tallied from zero, we have $g_0(z) \to g_l(z)$ and hold that $\big[[z^j]g_i(z)\big]=t_{i,j} \, \forall \, i,j$ which naturally means that 
\beq
  g_i(z)=\sum_{k=0}^{\infty} t_{i,k} z^k  \,\,\,\,\, \forall  \,\,\,\,\, 0 \leq i \leq l
\, .
\eeq
Further, it is held that one can create a reverse, partial sum matrix, through $g_l(z)*(1-z)^{-1}$ that is term-wise equivalent, so that the o.g.f. function created from the final difference sequence $\langle t_0,t_1,...,t_w \rangle_l$ as
\beq
  g_l(z)=\frac{t_{l,0}+t_{l,1}z+...+t_{l,w}z^w}{(1-z)^{l}}
\eeq
is the known o.g.f. that maps to the polynomial that defines the input sequence.
\end{lemma}
\begin{proof}
Given $(t_n)_{n \in \mathbb{N}}$, $t_n \in \mathbb{R}$, let us work on the terms of the first row, so that we have $g(z)(1-z)=(t_0z^0+t_1z^1+t_2z^2+t_3z^3+...)(1-z)$. This immediately gives us
{
\[ 
  \begin{tabular}{ccccccc}
       & $t_0z^0$ & +$t_1z^1$   & +$t_2z^2$      & +$t_3z^3$  & +$t_4z^4$ & ...  \\  
   $-$ &                  & \phantom{+}$t_0z^1$ & +$t_1z^2$     & +$t_2z^3$ & +$t_3z^4$ & ... \\
   \hline
       & $(t_0)z^0$ & +$(t_1-t_0)z^1$   & +$(t_2-t_1)z^2$ & +$(t_3-t_2)z^3$ & +$(t_4-t_3)z^4$ & ... \\
  \end{tabular}
\]
}
the very formula with which we defined the difference matrix at Definition~\ref{def:DifferenceSequenceMatrix}, which is to say $t_{i,j}=t_{i-1,j}-t_{i-1,j-1}$, and where \textit{off piste returns zero}.

Finally, the partial sum against generating series is assumed as common mathematical knowledge, and since a generating series defines one polynomial only, the two are equivalent, and uniquely so.
\end{proof}
\subsection{Deconvolution Matrices}
If we wish to generalise beyond the \textit{difference sequence key} we refer to the \textit{generating function key}, which is the numerator within the closed-form generating series. As regards deconvolution there is an interesting partial sum method for (not least) the Fibonacci series $z\! \cdot \!(1-z-z^2)^{-1}$ but beyond a small number of denominator terms it becomes impracticable, and we therefore require a more general deconvolution method.
\begin{lemma}
Let $\mathbb{R}[[z]]$ define the ring of formal power series in the variable $z$ over $\mathbb{R}$. Given a generating series $g(z)=a_0+a_1z+a_2z^2+...$ which has closed form 
\beq
  g(z)=\frac{b_0+b_1z+...+b_wz^w}{c_0+c_1z+...+c_uz^u}    
\, ,
\eeq
and assuming that $c_0 =1$, we can deduce $g(z)$ in open-form, and from the closed-form alone, via
\beq
\label{eq:MaclaurinShortCutMax}
  a_n=b_n-\sum_{k=1}^{\phi} c_k \! \cdot \! a_{n-k}
\, .
\eeq
Where $\phi=\min(u,n)$ and $u$ gives the highest power of the divisor, and where $n=0$ causes the sum to return zero.%
\footnote{Following decomposition norms one would expect $w<u$, and in fact earlier we defined midline sequences for polynomials at $w=u-1$ but this Lemma holds however $w$ and $u$ might relate, as long as there are some terms in the numerator and the denominator and $c_0 \neq 0$.}
\end{lemma}
\begin{example}
Given a closed form o.g.f. $(1+z^2)(1-2z-z^2+4z^3-z^4-2z^5+z^6)^{-1}$ the problem then is to discover the value of $a_7$.
{
\[
\begin{tabular}{cccccccc}
1&2&6&10&19&28&44 & ?\\
1&0&1&0  &0  &0  &0 & 0 \\
&&&&&& & $b_7$\\
\end{tabular}
\]
}
Thus we have
\beq
  a_{7}=b_{7}-\sum_{k=1}^{6} c_k a_{7-k}
\, ,
\eeq
which gives us
\beq
  a_{7}=b_{7}-(-2a_{6}-1a_{5}+4a_{4} -1a_{3}-2a_{2} +1a_{1})
\, .
\eeq
Which reduces to $a_{7}=0-(-88-28+76-10-12+2)$ and $a_{7}=60$.
\end{example}
A considerably faster method than applying the Maclaurin expansion on the closed-form generating series.
\begin{proof}
For deconvolution then one need look no further than a well known generating series formula for division. This says that if
\beq
  \sum_{n=0}^{\infty}b_n z^n \Big/ \sum_{n=0}^{\infty} c_n z^n=\sum_{n=0}^{\infty}a_n z^n
\, ,
\eeq
and that $g(z)=\sum_{n=0}^{\infty}a_n z^n$ with $\big[[z^n]g(z)\big]=a_n$ then
\beq
\label{eq:MaclaurinShortCut}
  a_n=\frac{1}{c_0}\bigg(b_n-\sum_{k=1}^n c_k \! \cdot \! a_{n-k} \bigg)
\, .
\eeq
With $c_0 \neq 0$, and where $n=0$ causes $\sum_{k=1}^n c_k \! \cdot \! a_{n-k}$ to return zero.

We note two things about $c_0$. Naturally it cannot be given by zero, but equally, if it is unity then it simplifies matters considerably. In combinatorics in particular it is not uncommon to assume that $c_0=1$. Now, whilst that is convenient $c_0$ can hold any value other than zero and the method will work equally well.

There are a number of ways that we can interpret this relation, however to conclude our proof the most expedient is to take $\sum^{\infty}_{n=0}b_nz^n$ as the numerator of the closed form o.g.f. and $\sum^{\infty}_{n=0}c_nz^n$ as the denominator. In this instance it so happens that the generating series are truncated by zeros.%
\footnote{For a non-truncated instance see Example~\ref{ex:DeconvolutionDifferenceOp}.} Indeed if we assume $c_0=1$ we obtain $a_n=b_n-\sum_{k=1}^{n} c_k a_{n-k}$ but since we needn't go beyond $u$ where $u$ gives the highest power of the divisor we obtain $a_n=b_n-\sum_{k=1}^{\phi} c_k\! \cdot \! a_{n-k}$ for $\phi=\min(u,n)$.
\end{proof}
\begin{example}
\label{ex:DeconvolutionDifferenceOp}
Thus we wish to reverse the difference operation i.e. convolution $(1-z)$, not with partial sums but with the pure deconvolution method. We begin with the result of a difference operation, so that if $g(z)=a_{0}+a_{1}z+a_{2}z^2+...$ then
\beq
  g(z)(1-z)=\big[ (a_0)+(a_1-a_0)z+(a_2-a_1)z^2+(a_3-a_2)z^3+(a_4-a_3)z^4+... \big]
\, .
\eeq
Now we wish to move \textit{from} the difference series in order to undo it, deconvolve it, to retrieve the input sequence. Hence the new terms are the $b_n$ values within our formula, where $b_0=a_0$, $b_1=a_1-a_0$ etc. Naturally we have $c_n$ for the divisor terms as before, where $c_0=1$ and $c_1=-1$. This significantly reduces the formula so that we now have
\beq
  a_{n}=b_{n}-(c_1\! \cdot \!a_{n-1})
\, ,
\eeq
which produces $a_0=a_0$, $a_1=(a_1-a_0)-(-a_0)$, $a_2=(a_2-a_1)-(-a_1)$ etc., as required.
\end{example}
\subsection{To Deconvolve Ramanujan's Closed-Form}
Within this section we would like to explore the deconvolution of $(1+z)^{-1}$ to provide clarity on $h(z)=\Lambda_s\! \cdot \!(1-z)^{-s-1}(1+z)^{-s-1}$ since we already know that the deconvolution $(1-z)^{-1}$ is a matrix of repeated partial sums on the numerator $\Lambda_s$.
\begin{proposition}
The deconvolution $(1+z)^{-1}$ produces an alternating partial sum of the original terms, where the first output term is identical to the first input term, and the second output term is the difference of the second and first input terms; with the former as minuend.
\end{proposition}
\begin{proof}
Given a generating series $g(z)=a_{0}+a_{1}z+a_{2}z^2+...$ and $g(z)(1+z)=b_0+b_1z+b_2z^2+b_3z^3+...$ let us apply a deconvolution per Equation~\ref{eq:MaclaurinShortCutMax}. Now, where we have $c_0=1$ and $c_1=1$, the formula reduces to
\beq
  a_{n}=b_{n}-a_{n-1}
\, ,
\eeq
and from this we can establish that $a_0=b_0$, $a_1=b_1-b_0$, $a_2=b_2-b_1+b_0$, $a_3=b_3-b_2+b_1-b_0$ etc.
\end{proof}
\begin{proposition}
Convolution against a unit convolvee by the pairs of operations $(1-z)^{-1}$, $(1+z)^{-1}$, through any number of iterations, will produce the pattern below.
{
\[
\begin{tabular}{>{\centering\arraybackslash}p{4mm}>{\centering\arraybackslash}p{4mm}>{\centering\arraybackslash}p{4mm}>{\centering\arraybackslash}p{4mm}>{\centering\arraybackslash}p{4mm}>{\centering\arraybackslash}p{4mm}>{\centering\arraybackslash}p{4mm}}
1&0&2&0&3&0&...\\
1&1&2&2&3&3&...\\
1&0&1&0&1&0&...\\
1&1&1&1&1&1&...\\
1&0&0&0&0&0&...\\
\end{tabular}
\]
}
Which is to say the same series that is produced from the application of $\frac{1}{(1-z)}$, for each pair of convolutions, but with zero terms inserted alternatively.
\end{proposition}
\begin{proof}
One need only perform the partial sum first, and examine the alternating partial sums second, in that order. They show us that odd numbered terms (since one tallies from zero) will cancel, and even valued terms will return $b_n$. What's more, even though these terms are dispersed by alternating zeros, it is clear that the partial sum will always return the appropriate figurate value.
\end{proof}
\subsection{Denouement}
Now for the conclusion of the proof, and to return to $g(z)$, the Eulerian rational functions, we observe that the cubes, from unity this time, can be constructed as 
\beq
  g(z)=\frac{1}{(1-z)^4}+\frac{4z}{(1-z)^4}+\frac{z^2}{(1-z)^4}
\, 
\eeq
which is the sum of multiple deconvolution matrices, carefully oriented. From our analysis it is clear that $h(z)=\frac{z+4z^3+z^5}{(1-z)^4(1+z)^4}$ too can be constructed in a similar manner
\beq
  h(z)=\frac{z}{(1-z)^4(1+z)^4}+\frac{4z^3}{(1-z)^4(1+z)^4}+\frac{z^5}{(1-z)^4(1+z)^4}
\, ,
\eeq
which is the sum of multiple deconvolution matrices, carefully oriented. And the overall pattern follows.
\begin{center}

\end{center}
  
  \vspace{12pt}
  \par
  \par  
  \textit{E-mail address}: \texttt{pjb15@student.london.ac.uk} \\

\end{document}